\documentclass[a4paper, oneside, 11pt]{article}
\usepackage[english]{babel}
\usepackage{amsfonts}
\usepackage{amssymb}
\usepackage{graphics}
\usepackage{amsrefs}
\usepackage{latexsym}
\begin{document}

\title{{\bf{\Large{A simple proof of a generalization of the Chu-Vandermonde identity}}}\footnote{{\it AMS (2000) subject classification}. Primary: 60G58. Secondary: 60G09.} }
\author{\textsc {Annalisa Cerquetti}\footnote{Corresponding author, SAPIENZA University of Rome, Via del Castro Laurenziano, 9, 00161 Rome, Italy. E-mail: {\tt annalisa.cerquetti@gmail.com}}\\
\it{\small Department of Methods and Models for Economics, Territory and Finance}\\
  \it{\small Sapienza University of Rome, Italy }}
\newtheorem{teo}{Theorem}
\date{\today}
\maketitle{}

\begin{abstract}
We provide a simple proof of a generalization of the multivariate Chu-Vandermo\-nde identity recently derived in Favaro et al. (2010a). Exploiting known results for rising factorials and fourth Lauricella polynomials we show resorting to Laplace-type integral representation of the fourth Lauricella function may be avoided.
\end{abstract}

\section {Introduction}

Favaro et al. (2010a) introduce the following generalized version of the multinomial theorem for rising factorials $(a)_n=a(a+1)\cdots(a+n-1)$ in terms of the fourth Lauricella function $F_D^{(r)}$ 
\begin{equation}
\label{multigen}
\sum_{n_1, \dots, n_r =0, \sum_j n_j=n}^{n} {n \choose n_1,\dots, n_r} \prod_{j=1}^r  w_j^{n_j} (a_j)_{n_j}$$$$= w_r^n (a)_n
 F_D^{(r-1)}\left[-n, a_1, \dots, a_{r-1}, a; 1-\frac{w_1}{w_r}, \dots, 1-\frac{w_{r-1}}{w_r} \right],
\end{equation} 
for $a_j$ and $w_j >0$, $j=1,\dots, r$, and $a=\sum_{j=1}^{r}a_j$, with application to species richness estimation within a Bayesian nonparametric framework outside the Gibbs {\it priors} class. The authors start by proving the binomial version of (\ref{multigen}) by repeated application of known properties of the Gauss hypergeometric {\it function}, and then obtain the multinomial version through a long proof resorting to integral representation of the fourth Lauricella {\it function}.
Here we show how to reduce both the proofs to few elementary steps of combinatorial calculus resorting to known properties of Gauss hypergeometric  and fourth Lauricella {\it polynomials} .\\

\section {The simple proof}
First we recall some elementary known facts about rising factorials and Gauss hypergeometric polynomials. See e.g. Srivastava and Manocha (1984). For $n \in \mathbb{N}$ and $x \neq 0, -1, -2,\dots,$
$$(x)_n=x(x+1)\cdots(x+n-1)=\frac{\Gamma(x +n)}{\Gamma(x)},$$
which implies for $0\leq k \leq n$ the binomial coefficient may be expressed as
\begin{equation}
\label{factor}
{n \choose k}= \frac{(-1)^k (-n)_k}{k!},
\end{equation}
as well as that
\begin{equation}
\label{gamrel}
\frac{\Gamma(n+1)}{\Gamma(n-k+1)}=(-1)^k (-n)_k
\end{equation}
and 
\begin{equation}
\label{ennek}
(x)_{n-k}=\frac{(-1)^k (x)_n}{(1-x-n)_k}.
\end{equation}
Now, for $a, b, c$ real or complex parameters, and $c \neq 0,-1,-2,\dots, $ the Gauss hypergeometric {\it function} (or Gauss hypergeometric series) in terms of rising factorials is given by the following infinite sum
\begin{equation}
\label{gauss}
_2 F_1(a, b, c; x)= \sum_{k=0}^\infty \frac{(a)_k (b)_k}{k! (c)_k}x^k
\end{equation}
and for $a$ (or $b$) a non positive integer reduces to the Gauss hypergeometric {\it polynomial}
\begin{equation}
\label{gaussn}
_2 F_1(-n, b, c; x)= \sum_{k=0}^n \frac{(-n)_k (b)_k}{k! (c)_k}x^k,
\end{equation}
which is known to satisfy the following relation (see e.g. Andrews et al. (1999) Eq. 2.3.14)
\begin{equation}
\label{relaz}
_2 F_1(-n, b, c; x)= \frac{(c-b)_n}{(c)_n} \hspace{1mm}_2 F_1 (-n, b, b+1 -n -c; 1-x).
\end{equation}
Notice that for $x=1$, both (\ref{gaussn}) and (\ref{relaz}) give the Gauss summation formula
$$
_2 F_1(-n, b, c; 1)= \frac{(c-b)_n}{(c)_n}.
$$
Hence by (\ref{factor}) and (\ref{ennek})
\begin{equation}
\label{prima}
\sum_{k=0}^n {n \choose k} (\alpha)_k w^k(\beta)_{n-k} z^{n-k} =(\beta)_n z^n \sum_{k=0}^n \frac{(-n)_k (\alpha)_k}{k! (1-\beta-n)_k}  \left(\frac w z\right)^k =
\end{equation}
and by (\ref{relaz}) and (\ref{gamrel})
\begin{equation}
= z^n(\alpha+\beta)_n \hspace{1mm} _2F_1 (-n, \alpha, \alpha+\beta, 1- \frac w z),
\end{equation} 
which proves the generalization of  the Chu-Vandermonde identity in Proposition 1. in Favaro et al. (2010a).\\\\

To prove the multivariate version (\ref{multigen}) it is enough and obvious to resort to the multivariate extension  of (\ref{relaz}) to the fourth Lauricella {\it polynomials} (Toscano, 1972, see Exton, 1976, pag. 216, Eq. 6.9.5)
\begin{equation}
\label{lauric}
F_D^{(r)}[-n, b_1,\dots, b_r; c; x_1, \dots, x_r]=$$
$$=\frac{(c - \sum_j b_j)_n}{(c)_n} F_D^{(r)}[-n , b_1, \dots, b_r; 1 + \sum_j b_j - n -c; 1-x_1, \dots, 1-x_r] 
\end{equation}
for  $$F_D^{(r)}= [-n, b_1, \dots, b_r; c; x_1, \dots, x_r]= \sum_{n_1,\dots, n_r=0}^{n, \sum_j n_j \leq n} \frac{(-n)_{\sum_j n_j}(b_1)_{n_1}\cdots (b_r)_{n_r}}{(c)_{\sum_j n_j}}\frac{x_1^{n_1}}{n_1!}\cdots\frac{x_r^{n_r}}{n_r!}.$$
Having at hand the generalizations of (\ref{factor}) 
\begin{equation}
\label{factormult}
\frac{n!}{n_1!\cdots n_r!}=\frac{(-1)^{\sum_{j=1}^{r-1} n_j} (-n)_{\sum_{j=1}^{r-1} n_j}}{n_1! \cdots n_{r-1}!}
\end{equation}
and (\ref{ennek})
\begin{equation}
\label{ennekmult}
(x)_{n -\sum_{j=1}^{r-1} n_j}=\frac{(-1)^{\sum_{j=1}^{r-1} n_j} (x)_{n}}{(1-x-n)_{\sum_{j=1}^{r-1} n_j}}
\end{equation}
the result (\ref{multigen}) (cfr. Proposition 2. in Favaro et al 2010a) easily follows by two steps of combinatorial calculus. In fact by (\ref{factormult}) and (\ref{ennekmult}), for $a=\sum_{j=1}^{r} a_j$
$$
\sum_{n_1, \dots, n_r =0, \sum_j n_j=n}^{n} {n \choose n_1,\dots, n_r} \prod_{j=1}^r  w_j^{n_j} (a_j)_{n_j}=$$
$$
=(a_r)_n w_r^n \sum_{n_1, \dots, n_r =0, \sum_j n_j=n}^{n} \frac{(-n)_{\sum_{j=1}^{r-1}  n_j}} {(1-a_r-n)_{\sum_{j=1}^{r-1} n_j}}\prod_{j=1}^{r-1} \frac{(a_j)_{n_j}}{n_j!} \left(\frac{w_j}{w_r}\right)^{n_j}= 
 $$
and by (\ref{lauric}) and (\ref{gamrel})
$$
= w_r^n \left(\sum_{j=1}^r a_j\right)_n F_D^{(r-1)} \left[-n, a_1, \dots, a_{r-1}; \sum_{j=1}^{r} a_j; 1-\frac{w_1}{w_r}, \dots, 1- \frac{w_{r-1}}{w_r}\right].
$$\\\
{\bf Remark} We notice that a lot of effort have been produced over the last few years around the multivariate Chu-Vandermonde identity in Bayesian nonparametric literature.  In Lijoi et al. (2007) and in Lijoi et al. (2008) two different proofs are provided for the multinomial theorem for rising factorial, despite the result is known in combinatorial theory, (see e.g. Carlson, 1977, pag 13, Th. 2.3-2). In  Favaro et al. (2010a) and in Favaro et al. (2010b) two different proofs are given to prove the generalization (1) despite as just shown it simply follows by two steps of combinatorial calculus.

\section*{References}
\newcommand{\bibu}{\item \hskip-1.0cm}
\begin{list}{\ }{\setlength\leftmargin{1.0cm}}

\bibu \textsc{Andrews, G.E., Askey, R. and Roy, R.} (1999) {\it Special functions}. In Encyclopedia of Mathematics and Its Applications. Cambridge University Press.

%\bibu \textsc{Cerquetti, A.} (2008) Generalized chinese restaurant construction of exchangeable Gibbs partitions and related results. {arXiv:0805.3853v1 [math.PR]}

\bibu \textsc{Carlson, B.C.} (1977) {\it Special functions of applied mathematics}. Academic Press, New York.

%\bibu \textsc{Cerquetti, A.} (2010) Bayesian nonparametric analysis of a species sampling model with finitely many types. {arXiv:1001.0245v1 [math.PR]}

\bibu \textsc{Exton, H.} (1976) {\it Multiple hypergeometric fuctions and applications}. Ellis Horwood, Chichester. 

\bibu\textsc{Favaro, S, Pr\"unster, I, and Walker, S. G.} (2010a) On a generalized Chu-Vandermonde Identity. {\it Methodol. Comput. Appl. Probab.} (to appear)

\bibu\textsc{Favaro, S, Pr\"unster, I, and Walker, S. G.} (2010b) On a class of random probability measures with general predictive structure. {\it Scand. J. Statist.} (to appear)

\bibu \textsc {Lijoi, A. Mena, R.H., Pr\"unster, I.} (2007) Bayesian nonparametric estimation of the probability of discovering a new species. {\it Biometrika} 94, 769--786

\bibu \textsc{Lijoi, A., Pr\"unster, I, Walker, S.G.} (2008) Bayesian nonparametric estimators derived from conditional Gibbs structures. {\it Ann. Appl. Probab.} 18:1519--1547.

%\bibu\textsc{Lauricella, G.} (1893) Sulle funzioni ipergeometriche a pi\`u  variabili. {\it Rend. Circ. Mat.} Palermo 7: 111--158. 

\bibu\textsc{Srivastava, H.M. and Manocha, H.L.} (1984) {\it A treatise on generating functions}. Wiley, New York.

\bibu\textsc{Toscano, L.} (1972) Sui polinomi ipergeometrici a pi\`u variabili del tipo FD di Lauricella. {\it Matematiche} (Catania), 27 , 219–250.

\end{list}
\end{document}